\documentclass[12pt]{article}
\newcommand{\be}{\begin{equation}}
\newcommand{\ee}{\end{equation}}
\newcommand{\bea}{\begin{eqnarray}}
\newcommand{\eea}{\end{eqnarray}}
\newcommand{\ba}{\begin{array}}
\newcommand{\ea}{\end{array}}

\newcommand{\bc}{\begin{center}}
\newcommand{\ec}{\end{center}}
\newcommand{\ben}{\begin{enumerate}}
\newcommand{\een}{\end{enumerate}}
\newcommand{\bfi}{\begin{figure}}
\newcommand{\efi}{\end{figure}}

\newcommand{\bq}{\begin{quote}}
\newcommand{\eq}{\end{quote}}
\newcommand{\bqu}{\begin{quotation}}
\newcommand{\equ}{\end{quotation}}
\newenvironment{emphit}{\begin{itemize}}{\end{itemize}}
\newcommand{\bemp}{\begin{emphit}}
\newcommand{\eemp}{\end{emphit}}

\newcommand{\bt}{\begin{tabular}}
\newcommand{\et}{\end{tabular}}

\newtheorem{myth}{Theorem}[section]

\newtheorem{mycor}{Corollary}[section]

\newtheorem{mydef}{Definition}[section]
\newtheorem{myrem}{Remark}[section]

\begin{document}
\date{}
\title{AN HOLOMORPHIC STUDY OF THE SMARANDACHE CONCEPT IN LOOPS
\footnote{2000 Mathematics Subject Classification. Primary 20NO5 ;
Secondary 08A05.}
\thanks{{\bf Keywords and Phrases : } holomorph of loops, Smarandache loops}}
\author{T\`em\'it\'op\'e Gb\'ol\'ah\`an Ja\'iy\'e\d ol\'a\thanks{On Doctorate Programme at
the University of Agriculture Abeokuta, Nigeria.}}
 \maketitle

\begin{abstract}
If two loops are isomorphic, then it is shown that their holomorphs
are also isomorphic. Conversely, it is shown that if their
holomorphs are isomorphic, then the loops are isotopic. It is shown
that a loop is a Smarandache loop if and only if its holomorph is a
Smarandache loop. This statement is also shown to be true for some
weak Smarandache loops(inverse property, weak inverse property) but
false for others(conjugacy closed, Bol, central, extra, Burn, A-,
homogeneous) except if their holomorphs are nuclear or central. A
necessary and sufficient condition for the Nuclear-holomorph of a
Smarandache Bol loop to be a Smarandache Bruck loop is shown.
Whence, it is found also to be a Smarandache Kikkawa loop if in
addition the loop is a Smarandache A-loop with a centrum holomorph.
Under this same necessary and sufficient condition, the
Central-holomorph of a Smarandache A-loop is shown to be a
Smarandache K-loop.
\end{abstract}

\section{Introduction}
The study of Smarandache loops was initiated by W.B. Vasantha
Kandasamy in 2002. In her book \cite{phd75}, she defined a
Smarandache loop (S-loop) as a loop with at least a subloop which
forms a subgroup under the binary operation of the loop. For more on
loops and their properties, readers should check \cite{phd3},
\cite{phd41},\cite{phd39}, \cite{phd49}, \cite{phd42} and
\cite{phd75}. In her book, she introduced over 75 Smarandache
concepts on loops. In her first paper \cite{phd83}, she introduced
Smarandache : left(right) alternative loops, Bol loops, Moufang
loops, and Bruck loops. But in this paper, Smarandache : inverse
property loops (IPL), weak inverse property loops (WIPL), G-loops,
conjugacy closed loops (CC-loop), central loops, extra loops,
A-loops, K-loops, Bruck loops, Kikkawa loops, Burn loops and
homogeneous loops will be introduced and studied relative to the
holomorphs of loops. Interestingly, Adeniran \cite{phd79} and
Robinson \cite{phd85}, Oyebo \cite{phdoyebo}, Chiboka and Solarin
\cite{phd80}, Bruck \cite{phd82}, Bruck and Paige \cite{phd40},
Robinson \cite{phd7}, Huthnance \cite{phd44} and Adeniran
\cite{phd79} have respectively studied the holomorphs of Bol loops,
central loops, conjugacy closed loops, inverse property loops,
A-loops, extra loops, weak inverse property loops and Bruck loops.\\

In this study, if two loops are isomorphic then it is shown that
their holomorphs are also isomorphic. Conversely, it is shown that
if their holomorphs are isomorphic, then the loops are isotopic.

It will be shown that a loop is a Smarandache loop if and only if
its holomorph is a Smarandache loop. This statement is also shown to
be true for some weak Smarandache loops(inverse property, weak
inverse property) but false for others(conjugacy closed, Bol,
central, extra, Burn, A-, homogeneous) except if their holomorphs
are nuclear or central. A necessary and sufficient condition for the
Nuclear-holomorph of a Smarandache Bol loop to be a Smarandache
Bruck loop is shown. Whence, it is found also to be a Smarandache
Kikkawa loop if in addition the loop is a Smarandache A-loop with a
centrum holomorph. Under this same necessary and sufficient
condition, the Central-holomorph of a Smarandache A-loop is shown to
be a Smarandache K-loop.

\section{Definitions and Notations}
Let $(L,\cdot )$ be a loop. Let $Aum(L,\cdot )$ be the automorphism
group of $(L,\cdot )$, and the set $H=(L,\cdot )\times Aum(L,\cdot
)$. If we define '$\circ$' on $H$ such that $(\alpha, x)\circ
(\beta, y)=(\alpha\beta, x\beta\cdot y)~\forall~(\alpha, x),(\beta,
y)\in H$, then $H(L,\cdot )=(H,\circ)$ is a loop as shown in Bruck
\cite{phd82} and is called the Holomorph of $(L,\cdot )$.

The nucleus of $(L,\cdot )$ is denoted by $N(L,\cdot )=N(L)$, its
centrum by $C(L,\cdot )=C(L)$ and center by $Z(L,\cdot )=N(L,\cdot
)\cap C(L,\cdot )=Z(L)$. For the meaning of these three sets,
readers should check earlier citations on loop theory.

If in $L$, $x^{-1}\cdot x\alpha\in N(L)$ or $x\alpha\cdot x^{-1}\in
N(L)~\forall~x\in L$ and $\alpha\in Aum(L,\cdot  )$, $(H,\circ)$ is
called a Nuclear-holomorph of $L$, if $x^{-1}\cdot x\alpha\in C(L)$
or $x\alpha\cdot x^{-1}\in C(L)~\forall~x\in L$ and $\alpha\in
Aum(L,\cdot  )$, $(H,\circ)$ is called a Centrum-holomorph of $L$
hence a Central-holomorph if $x^{-1}\cdot x\alpha\in Z(L)$ or
$x\alpha\cdot x^{-1}\in Z(L)~\forall~x\in L$ and $\alpha\in
Aum(L,\cdot  )$.

For the definitions of automorphic inverse property loop (AIPL),
anti-automorphic inverse property loop (AAIPL), weak inverse
property loop (WIPL), inverse property loop (IPL), Bol loop, Moufang
loop, central loop, extra loop, A-loop, conjugacy closed loop
(CC-loop) and G-loop, readers can check earlier references on loop
theory.

Here ; a K-loop is an A-loop with the AIP, a Bruck loop is a Bol
loop with the AIP, a Burn loop is Bol loop with the conjugacy closed
property, an homogeneous loop is an A-loop with the IP and a Kikkawa
loop is an A-loop with the IP and AIP.

\begin{mydef}
A loop is called a Smarandache inverse property loop (SIPL) if it
has at least a non-trivial subloop with the IP.

A loop is called a Smarandache weak inverse property loop (SWIPL) if
it has at least a non-trivial subloop with the WIP.

A loop is called a Smarandache G-loop (SG-loop) if it has at least a
non-trivial subloop that is a G-loop.

A loop is called a Smarandache CC-loop (SCCL) if it has at least a
non-trivial subloop that is a CC-loop.

A loop is called a Smarandache Bol-loop (SBL) if it has at least a
non-trivial subloop that is a Bol-loop.

A loop is called a Smarandache central-loop (SCL) if it has at least
a non-trivial subloop that is a central-loop.

A loop is called a Smarandache extra-loop (SEL) if it has at least a
non-trivial subloop that is a extra-loop.

A loop is called a Smarandache A-loop (SAL) if it has at least a
non-trivial subloop that is a A-loop.

A loop is called a Smarandache K-loop (SKL) if it has at least a
non-trivial subloop that is a K-loop.

A loop is called a Smarandache Moufang-loop (SML) if it has at least
a non-trivial subloop that is a Moufang-loop.

A loop is called a Smarandache Bruck-loop (SBRL) if it has at least
a non-trivial subloop that is a Bruck-loop.

A loop is called a Smarandache Kikkawa-loop (SKWL) if it has at
least a non-trivial subloop that is a Kikkawa-loop.

A loop is called a Smarandache Burn-loop (SBNL) if it has at least a
non-trivial subloop that is a Burn-loop.

A loop is called a Smarandache homogeneous-loop (SHL) if it has at
least a non-trivial subloop that is a homogeneous-loop.
\end{mydef}

\section{Main Results}
\subsection*{Holomorph of Smarandache Loops}
\begin{myth}\label{1:1}
Let $(L,\cdot )$ be a Smarandanche loop with subgroup $(S,\cdot )$.
The holomorph $H_S$ of $S$ is a group.
\end{myth}

\begin{myth}\label{1:2}
A loop is a Smarandache loop if and only if its holomorph is a
Smarandache loop.
\end{myth}
{\bf Proof}\\
Let $L$ be a Smarandache loop with subgroup $S$. By
Theorem~\ref{1:1}, $(H_S,\circ )$ is a group where $H_S=Aum(S,\cdot
)\times (S,\cdot )$. Clearly, $H_S\not \subset H(L,\cdot )$. So, let
us replace $Aum(S,\cdot )$ in $H_S$ by $A(S,\cdot )=\{\alpha\in
Aum(L,\cdot )~:~s\alpha\in S~\forall~s\in S\}$, the group of
Smarandache loop automorphisms on $S$ as defined in \cite{phd75}.
$A(S,\cdot )\le Aum(L,\cdot )$ hence, $H_S=A(S,\cdot )\times
(S,\cdot )$ remains a group. In fact, $(H_S,\circ )\subset
(H,\circ)$ and $(H_S,\circ )\leq (H,\circ)$. Thence, the holomorph
of a Smarandache loop is a Smarandache loop.

To prove the converse, recall that $H(L)=Aum(L)\times L$. If $H(L)$
is a Smarandache loop then $\exists~S_H\subset H(L)~\ni~S_H\leq
H(L)$. $S_H\subset H(L)\Rightarrow ~\exists~Bum(L)\subset Aum(L)$
and $B\subset L~\ni~S_H=Bum(L)\times B$. Let us choose
$Bum(L)=\{\alpha\in Aum(L)~:~b\alpha\in B~\forall~b\in B\}$, this is
the Smarandache loop automorphisms on $B$. So, $(S_H,\circ
)=(Bum(L)\times B,\circ )$ is expected to be a group.

Thus, $(\alpha ,x)\circ [(\beta ,y)\circ (\gamma ,z)]=[(\alpha
,x)\circ (\beta ,y)]\circ (\gamma ,z)~\forall~x,y,z\in
B,~\alpha,\beta,\gamma\in Bum(L)\Leftrightarrow x\beta\gamma\cdot
(y\gamma\cdot z)=(x\beta\gamma\cdot y\gamma )\cdot z\Leftrightarrow
x'\cdot (y'\cdot z)=(x'\cdot y')\cdot z~\forall~x',y',z\in B$. So,
$(B,\cdot )$ must be a group. Hence, $L$ is a Smarandache loop.

\begin{myrem}
It must be noted that if $Aum(L,\cdot )=A(S,\cdot )$, then $S$ is a
characteristic subloop.
\end{myrem}

\begin{myth}\label{1:3}
Let $L$ and $L'$ be loops. $L\cong L'$ implies $H(L)\cong H(L')$.
\end{myth}
{\bf Proof}\\
If $L\cong L'$ then $\exists$ a bijection $\alpha~:~L\to
L'~\ni~(\alpha ,\alpha ,\alpha )~:~L\to ~L'$ is an isotopism.
According to \cite{phd3}, if two loops are isotopic, then their
groups of autotopism are isomorphic. The automorphism group is one
of such since it is a form of autotopism. Thus ; $Aum(L)\cong
Aum(L')\Rightarrow H(L)=Aum(L)\times L\cong Aum(L')\times L'=H(L')$.

\begin{myth}\label{1:4}
Let $(L,\oplus)$ and $(L',\otimes )$ be loops. $H(L)\cong
H(L')\Leftrightarrow x\delta\otimes y\gamma =(x\beta\oplus
y)\delta~\forall~x,y\in L,~\beta\in Aum(L)$ and some
$\delta,\gamma\in Sym(L')$. Hence, $\gamma {\cal
L}_{e\delta}=\delta$, $\delta{\cal R}_{e\gamma}=\beta\delta$ where
$e$ is the identity element in $L$ and ${\cal L}_x$, ${\cal R}_x$
are respectively the left and right translations mappings of $x\in
L'$.
\end{myth}
{\bf Proof}\\
Let $H(L,\oplus )=(H,\circ )$ and $H(L',\otimes )=(H,\odot )$.
$H(L)\cong H(L')\Leftrightarrow ~\exists~\phi ~:~H(L)\to
H(L')~\ni~[(\alpha ,x)\circ (\beta ,y)]\phi =(\alpha,x)\phi\odot
(\beta ,y)\phi$. Define $(\alpha,x)\phi =(\psi^{-1}\alpha\psi
,x\psi^{-1}\alpha\psi )~\forall~(\alpha ,x)\in (H,\circ )$ and where $\psi~:~L\to L'$ is a bijection. \\
$H(L)\cong H(L')\Leftrightarrow (\alpha\beta ,x\beta\oplus y)\phi
=(\psi^{-1}\alpha\psi ,x\psi^{-1}\alpha\psi )\odot
(\psi^{-1}\beta\psi ,y\psi^{-1}\beta\psi )\Leftrightarrow
(\psi^{-1}\alpha\beta\psi ,(x\beta\oplus y)\psi^{-1}\alpha\beta\psi
)=(\psi^{-1}\alpha\beta\psi ,x\psi^{-1}\alpha\beta\psi\otimes
y\psi^{-1}\beta\psi )\Leftrightarrow (x\beta\oplus
y)\psi^{-1}\alpha\beta\psi=x\psi^{-1}\alpha\beta\psi\otimes
y\psi^{-1}\beta\psi\Leftrightarrow x\delta\otimes y\gamma
=(x\beta\oplus y)\delta$ where $\delta =\psi^{-1}\alpha\beta\psi$,
$\gamma =\psi^{-1}\beta\psi$.

Furthermore, $\gamma{\cal L}_{x\delta}=L_{x\beta}\delta$ and
$\delta{\cal R}_{y\gamma}=\beta R_y\delta~\forall~x,y\in L$. Thus,
with $x=y=e$, $\gamma{\cal L}_{e\delta}=\delta$ and $\delta{\cal
R}_{e\gamma}=\beta\delta$.

\begin{mycor}\label{1:5}
Let $L$ and $L'$ be loops. $H(L)\cong H(L')$ implies $L$ and $L'$
are isotopic under a triple of the form $(\delta ,I,\delta )$.
\end{mycor}
{\bf Proof}\\
In Theorem~\ref{1:4}, let $\beta =I$, then $\gamma =I$. The
conclusion follows immediately.

\begin{myrem}
By Theorem~\ref{1:3} and Corollary~\ref{1:5}, any two distinct
isomorphic loops are non-trivialy isotopic.
\end{myrem}

\begin{mycor}\label{1:6}
Let $L$ be a Smarandache loop. If $L$ is isomorphic to $L'$, then
$\{H(L), H(L')\}$ and $\{L,L'\}$ are both systems of isomorphic
Smarandache loops.
\end{mycor}
{\bf Proof}\\
This follows from Theorem~\ref{1:2}, Theorem~\ref{1:3},
Corollary~\ref{1:5} and the obvious fact that the Smarandache loop
property in loops is isomorphic invariant.

\begin{myrem}
The fact in Corollary~\ref{1:6} that $H(L)$ and $H(L')$ are
isomorphic Smarandache loops could be a clue to solve one of the
problems posed in \cite{phd83}. The problem required us to prove or
disprove that every Smarandache loop has a Smarandache loop
isomorph.
\end{myrem}

\subsection*{Smarandache Inverse Properties}
\begin{myth}\label{1:7}
Let $L$ be a loop with holomorph $H(L)$. $L$ is an IP-SIPL if and
only if $H(L)$ is an IP-SIPL.
\end{myth}
{\bf Proof}\\
In an IPL, every subloop is an IPL. So if $L$ is an IPL, then it is
an IP-SIPL. From \cite{phd82}, it can be stated that $L$ is an IPL
if and only if $H(L)$ is an IPL. Hence, $H(L)$ is an IP-SIPL.
Conversely assuming that $H(L)$ is an IP-SIPL and using the same
argument $L$ is an IP-SIPL.

\begin{myth}\label{1:8}
Let $L$ be a loop with holomorph $H(L)$. $L$ is a WIP-SWIPL if and
only if $H(L)$ is a WIP-SWIPL.
\end{myth}
{\bf Proof}\\
In a WIPL, every subloop is a WIPL. So if $L$ is a WIPL, then it is
a WIP-SWIPL. From \cite{phd44}, it can be stated that $L$ is a WIPL
if and only if $H(L)$ is a WIPL. Hence, $H(L)$ is a WIP-SWIPL.
Conversely assuming that $H(L)$ is a WIP-SWIPL and using the same
argument $L$ is a WIP-SWIPL.

\subsection*{Smarandache G-Loops}
\begin{myth}\label{1:9}
Every G-loop is a SG-loop.
\end{myth}
{\bf Proof}\\
As shown in [Lemma~2.2, \cite{phd45}], every subloop in a G-loop is
a G-loop. Hence, the claim follows.

\begin{mycor}\label{1:10}
CC-loops are SG-loops.
\end{mycor}
{\bf Proof}\\
In \cite{phd84}, CC-loops were shown to be G-loops. Hence, the
result follows by Theorem~\ref{1:9}.

\begin{myth}\label{1:11}
Let $G$ be a CC-loop with normal subloop $H$. $G/H$ is a SG-loop.
\end{myth}
{\bf Proof}\\
According to [Theorem~2.1,\cite{phd45}], $G/H$ is a G-loop. Hence,
by Theorem~\ref{1:9}, the result follows.

\subsection*{Smarandache Conjugacy closed Loops}
\begin{myth}\label{1:12}
Every SCCL is a SG-loop.
\end{myth}
{\bf Proof}\\
If a loop $L$ is a SCCL, then there exist a subloop $H$ of $L$ that
is a CC-loop. CC-loops are G-loops, hence, $H$ is a G-loop which
implies $L$ is a SG-loop.

\begin{myth}\label{1:13}
Every CC-loop is a SCCL.
\end{myth}
{\bf Proof}\\
By the definition of CC-loop in \cite{phd35}, \cite{phd47} and
\cite{phd78}, every subloop of a CC-loop is a CC-loop. Hence, the
conclusion follows.

\begin{myrem}
The fact in Corollary~\ref{1:10} that CC-loops are SG-loops can be
seen from Theorem~\ref{1:12} and Theorem~\ref{1:13}.
\end{myrem}

\begin{myth}\label{1:14}
Let $L$ be a loop with Nuclear-holomorph $H(L)$. $L$ is an
IP-CC-SIP-SCCL if and only if $H(L)$ is an IP-CC-SIP-SCCL.
\end{myth}
{\bf Proof}\\
If $L$ is an IP-CCL, then by Theorem~\ref{1:7}, $H(L)$ is an IP-SIPL
and hence by [Theorem~2.1, \cite{phd80}] and Theorem~\ref{1:13},
$H(L)$ is an IP-CC-SIP-SCCL. The converse is true by assuming that
$H(L)$ is an IP-CC-SIP-SCCL and using the same reasoning.

\subsection*{Smarandache : Bol loops, central loops, extra loops and Burn loops}
\begin{myth}\label{1:15}
Let $L$ be a loop with Nuclear-holomorph $H(L)$. $L$ is a Bol-SBL if
and only if $H(L)$ is a Bol-SBL.
\end{myth}
{\bf Proof}\\
If $L$ is a Bol-loop, then by \cite{phd85} and \cite{phd79}, $H(L)$
is a Bol-loop. According to [Theorem~6, \cite{phd83}], every
Bol-loop is a SBL. Hence, $H(L)$ is a Bol-SBL. The Converse is true
by using the same argument.

\begin{myth}\label{1:16}
Let $L$ be a loop with Nuclear-holomorph $H(L)$. $L$ is a
central-SCL if and only if $H(L)$ is a central-SCL.
\end{myth}
{\bf Proof}\\
If $L$ is a central-loop, then by \cite{phdoyebo}, $H(L)$ is a
central-loop. Every central-loop is a SCL. Hence, $H(L)$ is a
central-SCL. The Converse is true by using the same argument.

\begin{myth}\label{1:17}
Let $L$ be a loop with Nuclear-holomorph $H(L)$. $L$ is a extra-SEL
if and only if $H(L)$ is an extra-SEL.
\end{myth}
{\bf Proof}\\
If $L$ is a extra-loop, then by \cite{phd7}, $H(L)$ is a extra-loop.
Every extra-loop is a SEL. Hence, $H(L)$ is a extra-SEL. The
Converse is true by using the same argument.

\begin{mycor}\label{1:17.1}
Let $L$ be a loop with Nuclear-holomorph $H(L)$. $L$ is a
IP-Burn-SIP-SBNL if and only if $H(L)$ is an IP-Burn-SIP-SBNL.
\end{mycor}
{\bf Proof}\\
This follows by combining Theorem~\ref{1:14} and Theorem~\ref{1:15}.

\subsection*{Smarandache : A-loops, homogeneous loops}
\begin{myth}\label{1:18}
Every A-loop is a SAL.
\end{myth}
{\bf Proof}\\
According to [Theorem~2.2, \cite{phd40}], every subloop of an A-loop
is an A-loop. Hence, the conclusion  follows.

\begin{myth}\label{1:19}
Let $L$ be a loop with Central-holomorph $H(L)$. $L$ is an A-SAL if
and only if $H(L)$ is an A-SAL.
\end{myth}
{\bf Proof}\\
If $L$ is an A-loop, then by [Theorem~5.3, \cite{phd40}], $H(L)$ is
a A-loop. By Theorem~\ref{1:18}, every A-loop is a SAL. Hence,
$H(L)$ is an A-SAL. The Converse is true by using the same argument.

\begin{mycor}
Let $L$ be a loop with Central-holomorph $H(L)$. $L$ is an
homogeneous-SHL if and only if $H(L)$ is an homogeneous-SHL.
\end{mycor}
{\bf Proof}\\
This can be seen by combining Theorem~\ref{1:7} and
Theorem~\ref{1:19}.

\subsection*{Smarandache : K-loops, Bruck-loops and Kikkawa-loops}
\begin{myth}\label{1:19.1}
Let $(L,\cdot )$ be a loop with holomorph $H(L)$. $H(L)$ is an AIPL
if and only if $x\beta^{-1} J\cdot yJ=(x\cdot y\alpha^{-1}
)J~\forall~x,y\in L$ and $\alpha\beta =\beta\alpha~\forall~\alpha
,\beta\in Aum(L,\cdot )$. Hence, $xJ\cdot yJ=(z\cdot w)J~,~xJ\cdot
yJ=(x\cdot w)J~,~xJ\cdot yJ=(y\cdot w)J~,~xJ\cdot yJ=(z\cdot
x)J~,~xJ\cdot yJ=(z\cdot y)J~,~xJ\cdot yJ=(x\cdot y)J~,~xJ\cdot
yJ=(y\cdot x)J~\forall~x,y,z,w\in S$.
\end{myth}
{\bf Proof}\\
$H(L)$ is an AIPL $\Leftrightarrow~\forall~(\alpha ,x),(\beta ,y)\in
H(L)~,[(\alpha ,x)\circ (\beta ,y)]^{-1}=(\alpha ,x)^{-1}\circ
(\beta ,y)^{-1}\Leftrightarrow (\alpha\beta ,x\beta\cdot
y)^{-1}=(\alpha^{-1},(x\alpha^{-1})^{-1})\circ
(\beta^{-1},(y\beta^{-1})^{-1})\Leftrightarrow ((\alpha\beta )^{-1}
,[(x\beta\cdot y)(\alpha\beta
)^{-1}]^{-1})=(\alpha^{-1}\beta^{-1},(x\alpha^{-1})^{-1}\beta^{-1}\cdot
(y\beta^{-1})^{-1})\Leftrightarrow \alpha\beta
=\beta\alpha~\forall~\alpha ,\beta\in Aum(L,\cdot )$ and
$(x(\beta\alpha )^{-1})^{-1}\cdot
(y\beta^{-1})^{-1}=[x\alpha^{-1}\cdot y(\alpha\beta
)^{-1}]^{-1}\Leftrightarrow Aum(L,\cdot )$ is abelian and
$(x(\beta\alpha )^{-1})J\cdot y\beta^{-1} J=[x\alpha^{-1}\cdot
y(\alpha\beta )^{-1}]J\Leftrightarrow Aum(L,\cdot )$ is abelian and
$(x\alpha^{-1}\beta^{-1})J\cdot y\beta^{-1} J=[x\alpha^{-1}\cdot
y\beta^{-1}\alpha^{-1}]J\Leftrightarrow Aum(L,\cdot )$ is abelian
and $(x(\beta\alpha )^{-1})J\cdot y\beta^{-1} J=[x\alpha^{-1}\cdot
y(\alpha\beta )^{-1}]J\Leftrightarrow Aum(L,\cdot )$ is abelian and
$x'\beta^{-1}J\cdot y'J=(x'\cdot y'\alpha^{-1})J$ where
$x'=x\alpha^{-1},~y'=y\beta$.

What follows can be deduced from the last proof.

\begin{myth}\label{1:20}
Let $(L,\cdot )$ be a Bol-SBL with Nuclear-holomorph $H(L)$. $H(L)$
is a Bruck-SBRL if and only if $x\beta^{-1} J\cdot yJ=(x\cdot
y\alpha^{-1} )J~\forall~x,y\in L$ and $\alpha\beta
=\beta\alpha~\forall~\alpha ,\beta\in Aum(L,\cdot )$. Hence,
\begin{enumerate}
\item $L$ is a Moufang-SML and a Bruck-SBRL.
\item $H(L)$ is a Moufang-SML.
\item if $L$ is also an A-SAL with Centrum-holomorph $H(L)$ then $L$ is a Kikkawa-SKWL and so is
$H(L)$.
\end{enumerate}
\end{myth}
{\bf Proof}\\
By Theorem~\ref{1:15}, $H(L)$ is a Bol-SBL. So by
Theorem~\ref{1:19.1}, $H(L)$ is a Bruck-SBRL  $\Leftrightarrow
Aum(L,\cdot )$ is abelian and $x\beta^{-1} J\cdot yJ=(x\cdot
y\alpha^{-1} )J~\forall~x,y\in L$.
\begin{enumerate}
\item From Theorem~\ref{1:19.1}, $L$ is a Bruck-SBRL. From Theorem~\ref{1:19.1}, $L$ is an AAIPL, hence $L$
is a Moufang loop since it is a Bol-loop thus $L$ is a Moufang-SML.
\item $L$ is an AAIPL implies $H(L)$ is an AAIPL hence a Moufang loop. Thus, $H(L)$ is a
Moufang-SML.
\item If $L$ is also a A-SAL with Centrum-holomorph, then by Theorem~\ref{1:7}, $L$ and
$H(L)$ are both Kikkawa-Smarandache Kikkawa-loops.
\end{enumerate}

\begin{myth}\label{1:21}
Let $(L,\cdot )$ be a SAL with an A-subloop $S$ and
Central-holomorph $H(L)$. $H(L)$ is a SKL if and only if
$x\beta^{-1} J\cdot yJ=(x\cdot y\alpha^{-1} )J~\forall~x,y\in S$ and
$\alpha\beta =\beta\alpha~\forall~\alpha ,\beta\in A(S,\cdot )$.
Hence, $L$ is a SKL.
\end{myth}
{\bf Proof}\\
By Theorem~\ref{1:19}, $H(L)$ is a SAL with A-subloop $H_S=A(S,\cdot
)\times (S,\cdot )$. So $H(L)$ is a SKL if and only if $H_S$ is a
K-loop $\Leftrightarrow A(S,\cdot )$ is abelian and
$x\beta^{-1}J\cdot yJ=(x\cdot y\alpha^{-1})J~\forall~x,y\in
S~,\alpha,\beta\in A(S,\cdot )$ by Theorem~\ref{1:19.1}. Following
Theorem~\ref{1:19.1}, $S$ is an AIPL hence a K-loop which makes $L$
to be a SKL.

\paragraph{Address\\\\}

Mr. Jaiyeola Temitope Gbolahan,\\
Department of Mathematics, \\
Obafemi Awolowo University,\\
Ile Ife, Nigeria. \\
jaiyeolatemitope@yahoo.com \end{document}